\documentclass[a4paper]{article}
\usepackage{amsfonts}
\usepackage{amssymb}
\usepackage{latexsym}
\newtheorem{thm}{Theorem}[section]

\newtheorem{lem}[thm]{Lemma}
\newtheorem{rem}[thm]{Remark}
\newtheorem{dfn}[thm]{Definition}
\newtheorem{ex}[thm]{Example}

\newcommand{\proof}{\medskip \noindent {\bf Proof. \ \ }}
\newcommand{\qed}{\null\hfill $\Box\;\;$ \medskip}

\begin{document}

\parbox{1mm}

\begin{center}
{\bf {\sc \Large On vector-valued Dobrakov submeasures}}
\end{center}

\vskip 12pt

\begin{center}
{\bf Ondrej Hutn\'{\i}k}\footnote{{\it Mathematics Subject
Classification (2010):} 28B05, 28B15
\newline {\it Key words and phrases:} Non-additive set function,
Dobrakov submeasure, $L$-normed Banach lattice, vector-valued
measure, extension of a measure.
\newline {\it Acknowledgement.} This paper was supported by Grants
VEGA 2/0097/08 and VVGS 45/10-11.}
\end{center}

\vskip 24pt

\hspace{5mm}\parbox[t]{11cm}{\fontsize{9pt}{0.1in}\selectfont\noindent{\bf
Abstract.} Ivan Dobrakov has initiated a theory of non-additive
set functions defined on a ring of sets intended to be a
non-additive generalization of the theory of finite non-negative
countably additive measures. These set functions are now known as
the Dobrakov submeasures. In this paper we extend Dobrakov's
considerations to vector-valued submeasures defined on a ring of
sets. The extension of such submeasures in the sense of Drewnowski
is also given. } \vskip 24pt

\section{Introduction}

Non-additive set functions, as for example outer measures,
semi-variations of vector measures, appeared naturally earlier in
the classical measure theory concerning countable additive set
functions or more general finite additive set functions. A
systematic study of non-additive set function begins in the
fifties of the last century, cf.~\cite{Choquet}. Thence many
authors have investigated different kinds of non-additive set
functions, as submeasures~\cite{Drewnowski}, $t$-norms and
$t$-conorms~\cite{klement}, $k$-triangular set functions~\cite{AK}
and null-additive set functions~\cite{Pap}, fuzzy measures and
integrals~\cite{GMS, Mesiar} and many other types of set functions
and their properties. Specially, in different branches of
mathematics as potential theory, harmonic analysis, fractal
geometry, functional analysis, theory of nonlinear differential
equations, theory of difference equations and optimizations, etc.,
there are many types of non-additive set functions.

An interesting non-additive set function (as a generalization of a
notion of submeasure) was introduced by I.~Dobrakov.

\begin{dfn}\rm (Dobrakov, \cite{Dobrakov74})
Let $\mathcal{R}$ be a ring of subsets of a set $T \neq
\emptyset$. We say that a set function $\mu: \mathcal{R} \to [0,
+\infty)$ is a {\sl submeasure}, if it is
\begin{itemize}
\item[(1)] \label{monotone}\textsl{monotone}: if $A, B\in
\mathcal{R}$, such that $A \subset B$, then $\mu(A) \leq \mu(B)$;

\item[(2)] \label{continuous} \textsl{continuous at $\emptyset$}
(shortly \textsl{continuous}): for any sequence $(A_{n})_1^\infty$
of sets from $\mathcal{R}$, such that $A_{n} \searrow \emptyset$
(i.e., $A_n \supset A_{n+1}$ for each $n\in\mathbb{N}$ and
$\bigcap_{n\in \mathbb{N}} A_n = \emptyset$) there holds $\mu
(A_{n}) \to 0$ as $n\to\infty$;

\item[(3)] \label{subadditive}\textsl{subadditively continuous}:
for every $A \in \mathcal{R}$ and $\varepsilon>0$ there exists a
$\delta>0$, such that for every $B \in \mathcal{R}$ with
$\mu(B)<\delta$ there holds
\begin{enumerate}\item[(a)]
 $\mu(A \cup B) \leq \mu(A)+\varepsilon$, and
 \item[(b)]  $\mu(A) \leq \mu(A\setminus B)+\varepsilon$.
\end{enumerate}
\end{itemize}
\end{dfn}

Such a set function $\mu$ is now known as the {\sl Dobrakov
submeasure}. If the $\delta$ in condition~(3) is uniform with
respect to $A \in \mathcal{R}$, then we say that $\mu$ is a
\textsl{uniform Dobrakov submeasure}. Clearly, the definition of
Dobrakov submeasure provides a "non-additive generalization of the
theory of finite non-negative countably additive measures",
see~\cite{Dobrakov74}. If instead of~(3) we have $\mu(A \cup B)
\leq \mu(A) + \mu(B)$ for every $A,B \in \mathcal{R}$, or $\mu(A
\cup B) = \mu(A) + \mu(B)$ for every $A,B \in \mathcal{R}$ with $A
\cap B = \emptyset$, then we say that $\mu$ is a subadditive, or
an additive Dobrakov submeasure, respectively. Obviously,
subadditive, and particularly additive Dobrakov submeasures (i.e.,
countable additive measures) are uniform.

Note that there are two qualitative different types of continuity
of a set function $\mu$ in the definition. In literature, various
properties of continuity are added to the property~(1) in
Definition~\ref{monotone} when defining the notion of a submeasure
(and/or other generalizations, e.g. a semimeasure,
see~\cite{Dobrakov80}). There are also many papers where authors
consider various generalized settings (e.g.~\cite{Hal},
\cite{Hal2} and~\cite{Weber}). In paper~\cite{KlimkonSvistula}
authors considered the Darboux property of non-additive set
functions, in particular, the Dobrakov submeasure.
In~\cite{Riecan} and~\cite{Khare} we can find the (variant of)
Dobrakov submeasure in the context of fuzzy sets and systems.
In~\cite{Hal-Hut} some limit techniques to create new Dobrakov
submeasures from the old ones in the case when elements of the
ring $\mathcal{R}$ are subsets of the real line are developed. In
paper~\cite{Alyakin} Dobrakov submeasures with values in some
partially ordered semigroups are studied.

In this paper we extend the notion of a Dobrakov submeasure to set
functions with values in an $L$-normed Banach lattice (i.e., an
ordered space with a norm structure) and we investigate their
basic properties. Also, an extension theorem for the uniform
Dobrakov vector submeasures on a ring to a $\sigma$-ring is
discussed with respect to density in a topology induced by the
extended uniform Dobrakov vector submeasure. These results were
motivated by the work of Drewnowski~\cite{Drewnowski}.

\section{Preliminaries}\label{Section}

A \textsl{vector lattice} is a vector space equipped with a
lattice order relation, which is compatible with the linear
structure. A \textsl{Banach lattice} is defined to be a real
Banach space $\Xi$ which is also a vector lattice, such that the
norm $\|\cdot\|$ on $\Xi$ is monotone, i.e., $|x| \leq |y|$
implies $\|x\| \leq \|y\|$ for $x,y \in \Xi$, where for each $x
\in \Xi$ is $|x| = (x \vee 0) + (-x \vee 0)$ with $0$ being the
additive
identity on $\Xi$. The spaces $C(K)$, 
$L_p(\mu)$ for $1\leq p \leq +\infty$, 
and $c_0$ are important examples of Banach lattices.

A Banach lattice $\Xi$ is called an \textsl{abstract $L_1$-space}
(equivalently, an $L$-normed Banach lattice, or an $AL$-space) if
$\|x+y\| = \|x\|+\|y\|$ for all $x,y \geq 0$,
see~\cite{Birkhoff67} or~\cite{LinTza}. The spaces $L_1(\mu)$ and
$l_1$ are usual examples of $AL$-space.



An \textsl{order interval} $[x,y]$, where $x,y \in \Xi$, is the
set of all $z\in \Xi$, such that $x \leq z \leq y$. A subset $S
\subset \Xi$ is called \textsl{order bounded} if $S$ is contained
in some order interval of $\Xi$. A function $f: T \to \Xi$ is said
to be \textsl{order bounded} if its range is order bounded. If $f:
X \to Y$ and $Z \subset X$, then $f\mid_ Z$ is the restriction of
$f$ to $Z$.


In this paper $\Xi$ will represent an $AL$-space, and $\Lambda$
the positive cone of $\Xi$ (the set of all positive $(\geq)$
elements of $\Xi$). We also write $\overline{\Lambda} = \Lambda
\cup \{\lambda\}$, where $\lambda$ is such that $x < \lambda$ for
each $x \in \Xi$.

Let $\mathcal{R}$ be a collection of subsets of a non-void set $T$
which forms a ring under the operation $\triangle$ (symmetric
difference) and $\cap$ (intersection). As usual, a $\sigma$-ring
$\mathcal{S}$ is a collection of subsets of $T$ which is closed
under countable union and relative complementation. If
$\mathcal{A}, \mathcal{B} \subset \mathcal{R}$, then $$\mathcal{A}
\mathop{\cap}\limits^{\circ} \mathcal{B} = \{A \cap B; A \in
\mathcal{A}, B \in \mathcal{B}\}.$$ In the case
$\mathcal{A}=\{A\}$ we write $A\mathop{\cap}\limits^{\circ}
\mathcal{B}$ instead of $\{A\}\mathop{\cap}\limits^{\circ}
\mathcal{B}$. The operations $\mathop{\cup}\limits^{\circ},
\mathop{\triangle}\limits^{\circ}$ are defined similarly.

The following easy observations will be useful in the sequel of
this paper.

\begin{lem}\label{lemma1}
Let $\Lambda$ be the positive cone of an $AL$-space $\Xi$.
\begin{itemize} \item[(i)] If $\{f_{i}\} \subset
\Lambda$ is directed downward $(\geq)$ with $\inf_{i} f_{i} = f$,
where $f \in \Lambda$, then $\inf_{i} \|f_{i}\| = \|f\|$.
\item[(ii)] If $\{f_{i}\} \subset \Lambda$ is directed upward
$(\leq)$ with $\sup_{i} f_{i} = f$, where $f \in \Lambda$, then
$\sup_{i} \|f_{i}\| = \|f\|$.
\end{itemize}
\end{lem}

\proof Clearly, $\{f_{i} - f\} \in \Lambda$ is directed downward
$(\geq)$ with infimum $0$. Then according to results
in~\cite{Schaefer} (Ch.II, $\S$ 5.10 and Ch.II, $\S$ 1.7, $\S$ 2.4
and $\S$ 8.3) we have that $\lim_{i} \|f_{i} - f\| = 0.$ From it
follows that $\lim_{i} \|f_{i}\| = \|f\|$ and therefore $\inf_{i}
\|f_{i}\| = \|f\|$. The second item may be proved analogously.
\qed

Using these observations we immediately have the following

\begin{lem}\label{lemma2}
Let $\nu: \mathcal{M} \to \Lambda$ be a monotone set function,
where $\mathcal{M} \subset \mathcal{P}(T)$, $T \neq \emptyset$.
\begin{itemize}
\item[(i)] If $\mathcal{M}$ is closed with respect to finite
intersection, and $\inf \{\nu(A); E \subset A\in \mathcal{M}, E
\in T\} = a$, where $a\in \Lambda$, then $\inf \{\|\nu(A)\|; E
\subset A\in \mathcal{M}\} = \|a\|$. \item[(ii)] If $\mathcal{M}$
is closed with respect to finite union, and $\sup \{\nu(A); E
\supset A\in \mathcal{M}, E \in T\} = a$, where $a\in \Lambda$,
then $\sup \{\|\nu(A)\|; E \supset A\in \mathcal{M}\} = \|a\|$.
\end{itemize}
\end{lem}

\proof Let us prove the item~(i). It is obvious that the set
$P=\{\nu(A); E \subset A\in \mathcal{M}\}$ is a directed subset
$(\geq)$ of $\Lambda$, such that $\inf P = a$ exists in $\Lambda$.
From Lemma~\ref{lemma1}(i) we have that $\inf \{\|\nu(A)\|; E
\subset A\in \mathcal{M}\} = \|a\|$. The item~(ii) may be proved
similarly. \qed

\begin{dfn}\rm
The ordered pair $(\mathcal{R}, \Gamma)$, where $\mathcal{R}$ is a
ring and $\Gamma$ is a topology on $\mathcal{R}$, is called a
\textsl{topological ring of sets} if the ring operations $(A,B)
\to A\triangle B$ and $(A,B) \to A\cap B$ from $\mathcal{R} \times
\mathcal{R}$ (with the product topology) to $\mathcal{R}$ are
continuous.
\end{dfn}

The topology $\Gamma$ will be shortly called an $r$-topology on
$\mathcal{R}$. It it obvious that in a topological ring of sets
also the operations $(A,B) \to A\cup B$ and $(A,B) \to A \setminus
B$ are continuous. Recall that the notion of a topological ring of
sets is a generalization of spaces of measurable functions
introduced by Fr\'echet and Nikodym.

\begin{dfn}\rm
An $r$-topology $\Gamma$ on a ring $\mathcal{R}$ is said to be
\textsl{monotone}, or \textsl{Fr\'echet-Nikodym topology}
($FN$-topology, for short), if for each neighborhood $\mathcal{U}$
of $\emptyset$ there is a neighborhood $\mathcal{V}$ of
$\emptyset$, such that $\mathcal{V} \mathop{\cap}\limits^{\circ}
\mathcal{R} \subset \mathcal{U}$, i.e., such that $B \in
\mathcal{U}$ whenever $B \in \mathcal{R}$ and $B \subset A \in
\mathcal{V}$. A ring equipped with $FN$-topology is called an
\textsl{$FN$-ring}.
\end{dfn}

\begin{dfn}\rm
A base $\Omega$ at $\emptyset$ in $(\mathcal{R}, \Gamma)$ is
called a \textsl{normal base of neighborhoods of $\emptyset$} if
every $\mathcal{U} \in \Omega$ is a normal subclass of
$\mathcal{R}$ (i.e., $B\in \mathcal{U}$ provided $B\in
\mathcal{R}$ and $B \subset A$ for some $A\in \mathcal{U}$).
\end{dfn}

Now we introduce a notion of Dobrakov vector submeasure defined on
a ring $\mathcal{R}$ of subsets of a set $T\neq \emptyset$ with
values in an $AL$-space $\overline{\Lambda}$.

\begin{dfn}\rm\label{defvDobrak}
A set function $\mu: \mathcal{R} \to \overline{\Lambda}$ is called
a \textsl{Dobrakov vector submeasure}, briefly a
\textsl{$D$-submeasure}, if it is
\begin{itemize}
\item[(1)] \label{vmonotone}\textsl{monotone}: if $A, B\in
\mathcal{R}$, such that $A \subset B$, then $\mu(A) \leq \mu(B)$;

\item[(2)] \label{vcontinuous} \textsl{continuous}: for any
sequence $(A_{n})_1^\infty$ of sets from $\mathcal{R}$, such that
$A_{n} \searrow \emptyset$ there holds $\|\,\mu (A_{n})\| \to 0$
as $n\to\infty$;

\item[(3)] \label{vsubadditive}\textsl{subadditively continuous}
(s.c.): for every $A \in \mathcal{R}$ and $\varepsilon>0$ there
exists a $\delta>0$, such that for every $B \in \mathcal{R}$ with
$\|\,\mu(B)\|<\delta$ there holds
\begin{enumerate}\item[(a)]
 $\|\,\mu(A \cup B)\| \leq \|\,\mu(A)\|+\varepsilon$, and
 \item[(b)]  $\|\,\mu(A)\| \leq \|\,\mu(A\setminus B)\|+\varepsilon$.
\end{enumerate}
\end{itemize}
\end{dfn}

Note that the conditions~(3a) and~(3b) may be equivalently written
as the following sequence of inequalities
$$\|\,\mu(A)\|-\varepsilon \leq \|\,\mu(A \setminus B)\| \leq
\|\,\mu(A)\| \leq \|\,\mu(A \cup B)\| \leq
\|\,\mu(A)\|+\varepsilon.$$

Similarly as in the case of a Dobrakov submeasure, if the set
function $\mu$ has the property of \textsl{uniform subadditive
continuity}, shortly~(u.s.c.), then we say that $\mu$ is a
\textsl{uniform $D$-submeasure} ($D_{u}$-submeasure, for short).
If instead of~(3) we have $\|\,\mu(A \cup B)\| \leq \|\,\mu(A)\| +
\|\,\mu(B)\|$ for every $A,B \in \mathcal{R}$, or $\|\,\mu(A \cup
B)\| = \|\,\mu(A)\| + \|\,\mu(B)\|$ for every $A,B \in
\mathcal{R}$ with $A \cap B = \emptyset$, then we say that $\mu$
is a \textsl{subadditive $D$-submeasure} (shortly,
$D_{s}$-submeasure), or an \textsl{additive $D$-submeasure}
(shortly, $D_{a}$-submeasure), respectively.






\begin{ex}\rm\label{Cintegral}
Let $\mathcal{R}$ be a ring of subsets of $T \neq \emptyset$,
$T\in \mathcal{R}$, and $\mu: \mathcal{R} \to \Xi$ be a monotone
set function with $\mu(\emptyset) =0$ taking values in an
$AL$-space $\Xi$. Consider $f: T \to \mathbb{R}$ a non-negative
real function measurable with respect to $\mathcal{R}$ in the
sense $\{t \in T; f(t) > x\} \in \mathcal{R}$ for each $x\in
\mathbb{R}$. Analogously to~\cite{DHR} define the Choquet integral
of a function $f$ on a set $A$ with respect to $\mu$ by the
formula
$$\textrm{(C)} \int_{A} f\,d\mu = \int_{0}^{\infty} \mu(\{t \in A;
f(t) > x\})\,dx.$$ From the structural properties of set functions
defined by Choquet integral, see~\cite{Klir}, it is obvious that
if $\mu$ is a $D_s$-($D_a$-)submeasure, then the set function
$\nu_f: \mathcal{R} \to \Xi$ defined by $\nu_f(A) = \textrm{(C)}
\int_{A} f\,d\mu$ is also a $D_s$-($D_a$-)submeasure.

In this case the property~(s.c.) may be understood in the sense
that if two functions $f$ and $g$ differ on a set $A$ with measure
$\varepsilon$, then $\|\nu_f(A) - \nu_g(A)\| < \delta \cdot \tau$,
where $\tau = \sup_{t \in A}|f(t) - g(t)|$. Hence, we may estimate
errors in integration whenever we have some errors in inputs.
\end{ex}

\begin{rem}\rm
Observe that the integration technique developed in~\cite{Sipos,
Sipos2} may be extended to an $AL$-space $\Xi$ to obtain a
$\Xi$-{\sl valued \v Sipo\v s integral}. Recall that the \v Sipo\v
s integral is more general than the Choquet integral, but for
non-negative functions and fuzzy measures they coincide,
see~\cite{Pap}. The \v Sipo\v s integral is constructed as a limit
of nets. Such a case of Dobrakov net submeasures is investigated
in~\cite{Hal-Hut}. In particular, a $\Xi$-valued \v Sipo\v s
integral may also be considered as an example of Dobrakov vector
submeasure. Note that the \v Sipo\v s integral was successfully
used in prospect theory by Kahneman and Tversky, see~\cite{KT}. It
allows to describe how people make choices in situations where
they have to decide between alternatives involving risk.
\end{rem}

Concerning the notion of $D$-submeasure let us note that
the~(s.c.) in Definition~\ref{defvDobrak} may be replaced by the
following one.

\begin{lem}\label{lemmaequiv}
The set function $\mu: \mathcal{R} \to \overline{\Lambda}$ has
the~(s.c.) if and only if for $A,A_{n} \in \mathcal{R}$,
$n=1,2,\dots$, such that $\|\,\mu(A\triangle A_{n})\| \to 0$ holds
$\|\,\mu(A_{n})\| \to \|\,\mu(A)\|$ as $n\to\infty$.
\end{lem}

\proof Necessity: Suppose the contrary, i.e., let
$\|\,\mu(A_{n})\| \nrightarrow \|\,\mu(A)\|$ whenever
$\|\,\mu(A\triangle A_{n})\| \to 0$ for $A,A_{n} \in \mathcal{R}$,
$n=1,2,\dots$. Then we may assume that for some $\varepsilon >0$
either $\|\,\mu(A_{n})\| > \|\,\mu(A)\| + \varepsilon$ for each
$n\in \mathbb{N}$, or $\|\,\mu(A_{n})\| < \|\,\mu(A)\| -
\varepsilon$ for each $n\in \mathbb{N}$. In the first case we have
that $$\|\,\mu(A\cup (A\triangle A_{n}))\| \geq
\|\,\mu(A\triangle(A\triangle A_{n}))\|
> \|\,\mu(A)\| + \varepsilon,$$ which contradicts~(3a). Similarly
in the second case.

Sufficiency: Let $\|\,\mu(B_{n})\|\to 0$ as $n\to\infty$. Then
$$\|\,\mu(A \cup B_{n})\|= \|\,\mu(A \triangle (B_{n}\setminus A))\| \to \|\,\mu(A)\|,$$
and also $$\|\,\mu(A \setminus B_{n})\|= \|\,\mu(A \triangle
(B_{n}\cap A))\| \to \|\,\mu(A)\|$$ as $n\to\infty$. This
completes the proof. \qed

Lemma~\ref{lemmaequiv} may also be written as follows: a set
function $\mu: \mathcal{R} \to \overline{\Lambda}$ has the~(s.c.)
iff for each $A\in \mathcal{R}$ and each $\varepsilon >0$ there
exists a $\delta>0$, such that for each $C\in\mathcal{R}$ with
$\|\,\mu(A\triangle C)\,\|<\delta$ holds $\|\,\mu(C)\| -
\varepsilon < \|\,\mu(A)\| < \|\,\mu(C)\| + \varepsilon$.
Similarly we may prove that the property~(u.s.c.) is equivalent
with the following condition.

\begin{lem}\label{lemmaunisubadcont}
The set function $\mu: \mathcal{R} \to \overline{\Lambda}$ has
the~(u.s.c.) if and only if for $A_n$, $B_n \in \mathcal{R}$,
$n=1,2,\dots$, such that $\|\,\mu(A_n \triangle B_n)\,\| \to 0$
holds $\|\,\mu(A_n)\,\| - \|\,\mu(B_n)\,\| \to 0$ as $n\to\infty$.
\end{lem}

The property~(u.s.c.) says that for each $\varepsilon
>0$ there is a $\delta >0$, such that for all $A,B \in \mathcal{R}$
with $\|\,\mu(A \triangle B)\| < \delta$ holds $\|\,\mu(B)\| -
\varepsilon < \|\,\mu(A)\| < \|\,\mu(B)\| + \varepsilon$. For the
following definition see~\cite[Theorem 1]{Dobrakov80}.

\begin{dfn}\rm
A set function $\mu: \mathcal{R} \to \overline{\Lambda}$ is said
to have the \textsl{pseudometric generating property}, briefly
the~(p.g.p.), if for each $\varepsilon
>0$ there is a $\delta>0$, such that for every $A,B \in
\mathcal{R}$ with $\|\,\mu(A)\| \vee \|\,\mu(B)\| < \delta$ holds
$\|\,\mu(A \cup B)\| < \varepsilon$, where $a \vee b$, resp. $a
\wedge b$, means the maximum, resp. the minimum, of the real
numbers $a,b$.
\end{dfn}



\begin{ex}\rm
Consider the Choquet integral and $\nu_f(A) = \textrm{(C)}
\int_{A} f\,d\mu$. If $\|\nu_f(T)\| <+\infty$ and $\mu$ has
the~(p.g.p.), then $\nu_f$ has the~(p.g.p.) as well,
see~\cite{Li}.
\end{ex}

Clearly, the~(u.s.c.) implies the~(p.g.p.). The following theorem
rewritten in our setting is due to Dobrakov and Farkov\'a,
cf.~\cite[Lemma 3]{Dobrakov80}.

\begin{thm}\label{thmpgp}
Let $\mu: \mathcal{R} \to \overline{\Lambda}$ have the~(p.g.p.).
Then there is a sequence $(\delta_{k})_1^\infty$ of positive real
numbers with $\delta_k \searrow 0$, such that for any sequence
$(A_{k})_1^\infty$ of sets from $\mathcal{R}$ with
$\|\,\mu(A_{k})\| < \delta_{k}$ we have
$$\left\|\,\mu\left(\bigcup_{i=k+1}^{k+p} A_{i}\right)\right\|
< \delta_{k}$$ for each $k,p=1,2,\dots$.
\end{thm}

\proof Let $\mu$ have the~(p.g.p.). Then for $\varepsilon=1/2$
there exists a $\delta_1 \in (0,\frac{1}{2})$, such that for any
$A,B\in \mathcal{R}$  with $\|\,\mu(A)\| \vee \|\,\mu(B)\| <
\delta_1$ holds $\|\,\mu(A \cup B)\| < \frac{1}{2}$. For the above
$\delta_1$ there exists a $\delta_2 \in (0, \frac{1}{2^2} \wedge
\delta_1)$, such that for any $A,B\in \mathcal{R}$  with
$\|\,\mu(A)\| \vee \|\,\mu(B)\| < \delta_2$ we have $\|\,\mu(A
\cup B)\| < \delta_1$. Repeating this procedure we obtain a
sequence $(\delta_k)_1^\infty$, such that
$$0< \delta_{k+1} < \frac{1}{2^{k+1}} \wedge \delta_k, \quad
k=1,2,\dots.$$ If $\|\,\mu(A_k)\,\| < \delta_k$ for $k=1,2,\dots$,
then $$\left\|\,\mu\left(\bigcup_{i=k+1}^{k+p}
A_{i}\right)\right\| < \delta_{k}, \quad p=1,2,\dots.$$
\qed


\begin{dfn}\rm
A set function $\mu: \mathcal{R} \to \overline{\Lambda}$ is said
to be \textsl{exhaustive} on $\mathcal{R}$ if for each infinite
sequence $(A_{n})_1^\infty$ of pairwise disjoint sets from
$\mathcal{R}$ there holds $\|\,\mu(A_{n})\| \to 0$ as
$n\to\infty$.
\end{dfn}


\begin{dfn}\rm
Let $\mathcal{R}_{1}$ and $\mathcal{R}_{2}$ be two $\sigma$-rings,
such that $\mathcal{R}_{1} \subset \mathcal{R}_{2}$. If for every
$A \in \mathcal{R}_{2}$ there exists $B, C \in \mathcal{R}_{1}$,
such that $B \subset A \subset C$ and $\mu(C\setminus B) = 0$,
then $\mathcal{R}_{2}$ is called the \textsl{null-completion} of
$\mathcal{R}_{1}$.

We say that a $\sigma$-ring $\mathcal{S}$ is
\textsl{null-complete} with respect to $\mu$ if $B \subset A \in
\mathcal{S}$ and $\mu(A) = 0$, then $B \in \mathcal{S}$ and
$\mu(B) = 0$.
\end{dfn}

\section{Few elementary properties}

We begin with the following easy observations related to
$D_{s}$-submeasures on a ring.

\begin{thm}
Each $D_{s}$-submeasure $\mu$ on a ring $\mathcal{R}$ is
$\sigma$-subadditive, i.e.,
$$\left\|\,\mu\left(\bigcup_{n=1}^{\infty} A_{n}\right)\right\| \leq \sum_{n=1}^{\infty} \|\,\mu(A_{n})\|$$
for $A_{n} \in \mathcal{R}$, $n=1,2,\dots$.
\end{thm}

\proof Let $A_{n} \in \mathcal{R}$, $n=1,2,\dots$, such that
$\bigcup_{n=1}^{\infty} A_{n} = A \in \mathcal{R}$ and put $B_{n}
= A \setminus \bigcup_{i=1}^{n} A_{i}$, $n=1,2,\dots$. Then,
clearly, $B_{n} \in \mathcal{R}$, and $B_{n} \searrow \emptyset$.
Thus, $\|\,\mu(B_{n})\| \to 0$ as $n\to\infty$. Recall that if
$\mu$ is a $D_{s}$-submeasure on $\mathcal{R}$, then
$$\left\|\,\mu\left(\bigcup_{i=1}^{n}A_{i}\right)\right\| \leq \sum_{i=1}^{n} \|\,\mu(A_{i})\|$$
for every finite sequence $(A_{i})_1^n$ of arbitrary sets from
$\mathcal{R}$. Since $A \subset B_{n} \cup \bigcup_{i=1}^{n}
A_{i}$ for every $n\in \mathbb{N}$, then we get
{\setlength\arraycolsep{2pt}
\begin{eqnarray*}
\|\,\mu(A)\| & \leq & \left\|\,\mu\left(\bigcup_{i=1}^{n} B_{n}
\cup A_{i}\right)\right\| \leq \sum_{i=1}^{n} \|\,\mu(B_{n} \cup
A_{i})\|
\\ & \leq & \|\,\mu(B_{n})\| + \sum_{i=1}^{n} \|\,\mu(A_{i})\|.
\end{eqnarray*}}From it follows $$\|\,\mu(A)\| \leq \lim_{n\to \infty} \|\,\mu(B_{n})\| +
\sum_{i=1}^{\infty} \|\,\mu(A_{i})\| = \sum_{i=1}^{\infty}
\|\,\mu(A_{i})\|.$$ Hence the result. \qed

\begin{thm}\label{thmRsigma}
Let $\mu$ be a $D$-submeasure on $\mathcal{R}$ and
$(A_{n})_1^\infty$ be a sequence of sets from $\mathcal{R}$, such
that $A_{n} \nearrow (\searrow) A$, $A \in \mathcal{R}$. Then
$$\|\,\mu(A)\| = \|\,\mu(\lim_{n\to \infty} A_{n})\| = \lim_{n\to
\infty} \|\,\mu(A_{n})\|.$$
\end{thm}

\proof Suppose that $A_{n}\nearrow A$. Then $A \triangle A_{n} = A
\setminus A_{n}$ and obviously $A \setminus A_{n} \searrow
\emptyset$. From continuity of $\mu$ we have that $\|\,\mu(A
\setminus A_{n})\| \to 0$ as $n\to\infty$, and therefore
$\|\,\mu(A \triangle A_{n})\| \to 0$ as $n\to\infty$. Using
Lemma~\ref{lemmaequiv} we immediately get $\|\,\mu(A_{n})\| \to
\|\,\mu(A)\|$, i.e.,
$$\lim_{n\to \infty} \|\,\mu(A_{n})\| = \|\,\mu(A)\| =
\|\,\mu(\lim_{n\to \infty} A_{n})\|.$$ Analogously we may prove
the result for $A_{n} \searrow A$. \qed

\begin{thm}\label{thmex}
A $D$-submeasure $\mu$ is exhaustive on a ring $\mathcal{R}$ if
and only if every monotone sequence $(A_{n})_1^\infty$ of sets
from $\mathcal{R}$ is $\mu$-Cauchy, i.e., $$\|\,\mu(A_{n}
\triangle A_{m})\| \to 0\,\,\textrm{whenever}\,\,n \wedge m \to
\infty.$$
\end{thm}

\proof Necessity: Suppose the contrary, i.e., let
$(A_{n})_1^\infty$ be a monotone sequence of sets from
$\mathcal{R}$ which is not $\mu$-Cauchy. Without loss of
generality let us assume that the sequence $(A_n)_1^\infty$ is
increasing. Then there exists a positive integer $N$ and (an
infinite number of) $n_1, n_2, \dots$, where $n_i > N$,
$i=1,2,\dots$, such that $\|\,\mu(A_{n_j} \triangle A_{n_k})\|
\geq \varepsilon$ for $j\neq k$. We set
$$P_{n_{k}} = A_{n_{k+1}}\triangle A_{n_{k}} =
A_{n_{k+1}}\setminus A_{n_{k}}.$$ Clearly, $P_{n_{k}} \cap
P_{n_{k+1}} = \emptyset$ for $k=1,2,\dots$. Now,
$(P_{n_{k}})_1^\infty$ is a disjoint sequence of sets from
$\mathcal{R}$, such that $\|\,\mu(P_{n_{k}})\| \geq \varepsilon$
for $k=1,2,\dots$. This contradicts the fact that $\mu$ is
exhaustive.

Sufficiency: Let $(A_{n})_1^\infty$ be a disjoint sequence of sets
from $\mathcal{R}$ and put $B_n = \bigcup_{k=1}^{n} A_k$. If
$\|\mu(A_n)\|\nrightarrow 0$ as $n\to\infty$, there exists an
$\varepsilon
> 0$ and an increasing sequence $(n_k)_1^\infty$ of natural numbers, such that
$\|\,\mu(A_{n_k})\,\| > \varepsilon$ for $k = 1, 2, \dots$. Then
$\|\,\mu(B_{n_k})\,\| \geq \|\,\mu(A_{n_k})\,\| > \varepsilon$ for
$k = 1, 2, \dots$, which contradicts the fact that
$\|\,\mu(B_{n_k})\,\|$ is Cauchy. \qed


The following result shows that the situation from
Theorem~\ref{thmex} is different when considering a $D$-submeasure
on a $\sigma$-ring.

\begin{thm}\label{thmDsubsigmaring}
Each $D$-submeasure $\mu: \mathcal{S} \to \overline{\Lambda}$ on a
$\sigma$-ring $\mathcal{S}$ is exhaustive. 
\end{thm}

\proof Let $(A_{n})_1^\infty$ be a disjoint sequence of sets from
$\mathcal{S}$ and put $B_{n} = \bigcup_{k=n}^{\infty} A_{k}$. Then
$B_{n} \searrow \emptyset$, and from continuity of $\mu$ we have
$\|\,\mu(B_{n})\| \to 0$ as $n\to\infty$. Since $\mu(A_{n}) \leq
\mu(B_{n})$ for every $n\in \mathbb{N}$, then it follows that
$\|\,\mu(A_{n})\|
\to 0$ as $n\to\infty$. Thus $\mu$ is exhaustive on $\mathcal{S}$. 
\qed

\begin{thm}
Let $\mu: \mathcal{R} \to \overline{\Lambda}$ be an order bounded
$D_{u}$-submeasure on a ring $\mathcal{R}$. Then the class
$\mathcal{T}$ of all $\mathcal{U}_{\varepsilon}$
$(0<\varepsilon)$, where $\mathcal{U}_{\varepsilon} = \{A\in
\mathcal{R}; \|\,\mu(A)\| \leq \varepsilon\}$, forms a normal base
of neighborhoods at $\emptyset$ for an $FN$-topology.
\end{thm}

\proof It is easy to see that $\mathcal{T}$ is a filter base
satisfying the following conditions
\begin{itemize}
\item[(1)] for each $\mathcal{U} \in \mathcal{T}$ there exists
$\mathcal{V} \in \mathcal{T}$, such that $\mathcal{V}
\mathop{\triangle}\limits^{\circ} \mathcal{V} \subset
\mathcal{U}$; \item[(2)] for each $\mathcal{U} \in \mathcal{T}$
there exists $\mathcal{V} \in \mathcal{T}$, such that $\mathcal{V}
\mathop{\cap}\limits^{\circ} \mathcal{V} \subset \mathcal{U}$;
\item[(3)] for each $A \in \mathcal{R}$ and $\mathcal{U} \in
\mathcal{T}$ there exists $\mathcal{V} \in \mathcal{T}$, such that
$A \mathop{\cap}\limits^{\circ} \mathcal{V} \subset \mathcal{U}$.
\end{itemize} From the general theory of topological rings~\cite{Bourbaki} and
according to~\cite[$\S 1$]{Drewnowski} these three conditions are
necessary and sufficient that a filter base $\mathcal{T}$ of
neighborhoods of $\emptyset$ determines an $r$-topology on
$\mathcal{R}$. It is clear, that this topology is an
$FN$-topology. Moreover, the filter base $\mathcal{T}$ has the
following properties
\begin{itemize}
\item[(4)] each class $\mathcal{U}\in \mathcal{T}$ is normal in
$\mathcal{R}$, and \item[(5)] for each $\mathcal{U} \in
\mathcal{T}$ there exists $\mathcal{V} \in \mathcal{T}$, such that
$\mathcal{V} \mathop{\cup}\limits^{\circ} \mathcal{V} \subset
\mathcal{U}$.
\end{itemize}
Then according to~\cite[p. 142]{Vladimirov} $\mathcal{T}$ is a
normal base of neighborhoods of $\emptyset$ for an $FN$-topology
generated (or determined) by $\mu$ on $\mathcal{R}$. \qed

\begin{rem}\rm\label{remark3.6}
The $FN$-topology generated by $\mu$ on $\mathcal{R}$ is denoted
by $\Gamma(\mu)$. Since the concept of (s.c.) of $\mu$ is linked
with absolute continuity, in fact, only the continuity of $\mu$
and the condition (a.c.)
$$\|\,\mu(A_{n})\| + \|\,\mu(B_{n})\| \to 0 \,\, \Rightarrow \,\, \|\,\mu(A_{n} \cup B_{n}) \| \to
0$$ as $n\to\infty$ are needed for $\Gamma(\mu)$ to be an
$FN$-topology, see~\cite{DrewnowskiCM}. Clearly,
$D_{u}$-submeasures satisfy this condition. On the other hand,
$D$-submeasures do not satisfy the~(a.c.) in general.
\end{rem}

To prove the next theorem we first recall two Drewnowski's results
from~\cite{Drewnowski}.

\begin{lem}\label{lemmadrew1}
If $(\mathcal{R}, \Gamma)$ is a topological ring of sets and
$\mathcal{P}$ is a subring of the ring $\mathcal{R}$, then
$\overline{\mathcal{P}}^{\Gamma}$ is a subring of $\mathcal{R}$,
where $\overline{\mathcal{P}}$ denotes the closure of
$\mathcal{P}$ in $(\mathcal{R}, \Gamma)$.
\end{lem}

\begin{lem}\label{lemmadrew2}
If $(\mathcal{R}, \Gamma)$ is a topological ring of sets and
$\Omega$ is a base of (the filter of all) neighborhoods of
$\emptyset$ in $\mathcal{R}$, then for each $A \in \mathcal{R}$,
$A \triangle \Omega = \{A \triangle \mathcal{U}; \mathcal{U}\in
\Omega\}$ is a base of (the filter of all) neighborhoods of $A$ in
$\mathcal{R}$.
\end{lem}

\begin{thm}\label{thmdense}
Let $\sigma(\mathcal{R})$ be a $\sigma$-ring generated by a ring
$\mathcal{R}$ and let $\mu$ be an order bounded $D_{u}$-submeasure
on $\sigma(\mathcal{R})$. Then $\mathcal{R}$ is dense in
$(\sigma(\mathcal{R}), \Gamma(\mu))$.
\end{thm}

\proof Denote by $\overline{\mathcal{R}} =
\overline{\mathcal{R}}^{\Gamma(\mu)}$. According to
Lemma~\ref{lemmadrew1} we have that $\overline{\mathcal{R}}$ is a
subring of $\sigma(\mathcal{R})$.

Let $(A_{n})_1^\infty$ be a disjoint sequence of sets from
$\overline{\mathcal{R}}$, such that $\bigcup_{n=1}^{\infty} A_{n}=
A$. Then obviously, $$B_{n} = \bigcup_{k=1}^{n} A_{k} \in
\overline{\mathcal{R}}, \quad \textrm{for every \,} n\in
\mathbb{N}.$$ Put $$C_{n} = A \triangle B_{n} = A \triangle
\left(\bigcup_{k=1}^{n} A_{k}\right) = \bigcup_{k=n+1}^{\infty}
A_{k}.$$ Clearly, $C_{n} \searrow \emptyset$. Let $\varepsilon >0$
and $$\mathcal{V} = \left\{E\in \sigma(\mathcal{R});\,\,
\|\,\mu(E)\| \leq \frac{\varepsilon}{2}\right\}$$ be a
neighborhood of $\emptyset$ in $\sigma(\mathcal{R})$. Then for
each $n \in \mathbb{N}$ the neighborhood $B_{n} \triangle
\mathcal{V}$ of $B_{n}$ contains an element $E_{n} = B_{n}
\triangle V_{n} \in \mathcal{R}$, where $V_{n} \in \mathcal{V}$,
and also $$\|\,\mu(A \triangle E_{n})\| = \|\,\mu(C_{n} \triangle
V_{n})\| \leq \|\,\mu(C_{n} \cup V_{n})\|.$$ From continuity of
$\mu$ we have that $\|\,\mu(C_{n})\| \to 0$ as $n\to\infty$, and
therefore
$$\|\,\mu(C_{n} \cup V_{n})\| \leq \|\,\mu(V_{n})\| +
\frac{\varepsilon}{2},$$ which is possible by the~(u.s.c.) of
$\mu$. Since $V_{n} \in \mathcal{V}$, then $\|\,\mu(V_{n})\| \leq
\frac{\varepsilon}{2}$ for every $n=1,2,\dots$, and therefore
$$\|\,\mu(A \triangle E_{n})\| \leq \|\,\mu(C_{n} \cup
V_{n})\| \leq \|\,\mu(V_{n})\| + \frac{\varepsilon}{2} \leq
\frac{\varepsilon}{2} + \frac{\varepsilon}{2} = \varepsilon.$$
Since $A \triangle E_{n} \in \sigma(\mathcal{R})$ for all $n\in
\mathbb{N}$, then $A \triangle E_{n} \in
\mathcal{U}_{\varepsilon}$, where
$$\mathcal{U}_{\varepsilon} = \{F\in \sigma(\mathcal{R}); \,\,
\|\,\mu(F)\| \leq \varepsilon\}
$$ is a neighborhood of $\emptyset$ in $\sigma(\mathcal{R})$. Accordingly,
$E_{n} = A \triangle (A \triangle E_{n}) \in A \triangle
\mathcal{U}_{\varepsilon}$. Therefore each neighborhood of $A$
contains an element of $\mathcal{R}$ (according to
Lemma~\ref{lemmadrew2}). Hence $A\in \overline{\mathcal{R}}$, and
therefore $\overline{\mathcal{R}}$ is a $\sigma$-ring. Thus,
$\overline{\mathcal{R}}=\sigma(\mathcal{R})$. This completes the
proof. \qed

\section{Extension of D-submeasure}

In measure theory, an essential concept is the extension of the
notion of a measure (or, a submeasure) on one class of sets to a
notion of measure (or, a submeasure) on a larger class of sets.
For instance, in~\cite{Dobrakov84} Dobrakov showed the following
extension of a (Dobrakov) submeasure from a ring to a generated
$\sigma$-ring: \emph{An additive, subadditive or uniform
(Dobrakov) submeasure $\mu: \mathcal{R} \to [0,+\infty)$ has a
unique extension $\mu: \sigma(\mathcal{R}) \to [0,+\infty)$ of the
same type if and only if $\mu$ is exhaustive.} In this section we
study the possibility of an extension for a $D_{u}$-submeasure
defined on a ring $\mathcal{R}$ to a $\sigma$-ring
$\mathcal{R}_{0}$ in the sense that $\mathcal{R}$ is dense in
$\mathcal{R}_{0}$ with respect to a topology induced by the
extended $D_{u}$-submeasure.

Let $\mathcal{R}$ be a ring of subsets of $T \neq \emptyset$. Then
$$\mathcal{R}_{\sigma} = \{A;\, \textrm{there are \,} A_{n}\in \mathcal{R}, n=1,2,\dots, \,\,
\textrm{such that \,} A_{n}\nearrow A \}$$ denotes the standard
class of limits of increasing sequences of sets of $\mathcal{R}$.
It is clear that $\mathcal{R}_{\sigma}$ is closed with respect to
countable unions and finite intersections. Also, if $A \in
\mathcal{R}_{\sigma}$ and $B \in \mathcal{R}$, then $A\setminus B
\in \mathcal{R}_{\sigma}$.


Let $\mu: \mathcal{R} \to \overline{\Lambda}$ be an order bounded
exhaustive $D_{u}$-submeasure on a ring $\mathcal{R}$ and for each
$A \in \mathcal{R}_{\sigma}$ define the set function $\hat{\mu}:
\mathcal{R}_{\sigma} \to \overline{\Lambda}$ as follows
\begin{equation}\label{hatmu}
\hat{\mu}(A) = \sup\{\mu(B); B \subset A, B \in \mathcal{R}\}.
\end{equation} If $(C_{n})_1^\infty$ is a sequence of sets from $\mathcal{R}$, such that
$A=\bigcup_{n=1}^{\infty} C_{n}$, then there exists a sequence
$(B_{n})_1^\infty$ of sets from $\mathcal{R}$ with $B_{1} \subset
B_{2} \subset \dots$, such that
$$B_{n}=\bigcup_{i=1}^{n}C_{i} \quad \textrm{and} \quad
\bigcup_{n=1}^{\infty} B_{n} = \bigcup_{n=1}^{\infty} C_{n} = A.$$
From Lemma~\ref{lemma2}(ii) it follows that
$$\|\,\hat{\mu}(A)\| = \sup\{\|\,\mu(B)\|; B \subset A, B\in \mathcal{R}\}.$$
Then it is obvious that $$\|\,\hat{\mu}(A)\| =
\sup\{\|\,\mu(B_{n})\|; B_{n} \subset A, B_{n} \nearrow A,
B_{n}\in \mathcal{R}\},$$ which results
\begin{equation}\label{eqhatmu}
\|\,\mu(B_{n})\| \to \|\,\hat{\mu}(A)\| \quad \textrm{as\,\,} n\to
\infty.
\end{equation}

\begin{thm}\label{thmhatmu}
Let $\mu: \mathcal{R} \to \overline{\Lambda}$ be an order bounded
exhaustive $D_{u}$-submeasure on a ring $\mathcal{R}$ and
$\hat{\mu}: \mathcal{R}_{\sigma} \to \overline{\Lambda}$ be
defined as in~(\ref{hatmu}). Then $\hat{\mu}$ has the following
properties:
\begin{itemize}
\item[(a)] $\hat{\mu} \mid_ \mathcal{R} = \mu$, $\hat{\mu}$ is
monotone; \item[(b)] $\hat{\mu}$ is exhaustive on
$\mathcal{R}_{\sigma}$; \item[(c)] if $A_{n} \in \mathcal{R}$,
$n=1,2,\dots$, such that $A_{n} \nearrow A$, then
$\|\,\hat{\mu}(A\setminus A_{n})\| \to 0$ as $n\to\infty$;
\item[(d)] $\hat{\mu}$ has the~(u.s.c.) on $\mathcal{R}_{\sigma}$;
\item[(e)] $\hat{\mu}$ is continuous on $\mathcal{R}_{\sigma}$.
\end{itemize}
\end{thm}

\proof The item~(a) is obvious.

(b) Let $(A_{n})_1^\infty$ be a disjoint sequence of sets from
$\mathcal{R}_{\sigma}$. We have that
$$\|\,\hat{\mu}(A_{n})\| = \sup\{\|\,\mu(C)\|;\, C \subset A_{n}, C \in
\mathcal{R}\}.$$ Let $\varepsilon >0$ be chosen arbitrarily. Then
there exists $B_{n} \in \mathcal{R}$, such that $B_{n} \subset
A_{n}$ and
$$\|\,\hat{\mu}(A_{n})\| < \|\,\mu(B_{n})\| +
\frac{\varepsilon}{2^{n}}, \quad n=1,2,\dots .$$ Since
$(A_{n})_1^{\infty}$ is a disjoint sequence, then
$(B_{n})_1^{\infty}$ is disjoint as well. Also, $\mu$ is
exhaustive on $\mathcal{R}$, i.e., $\|\,\mu(B_{n})\| \to 0$ as
$n\to\infty$. Consequently, $\|\,\hat{\mu}(A_{n})\| \to 0$ as
$n\to\infty$ and thus, $\hat{\mu}$ is exhaustive on
$\mathcal{R}_{\sigma}$.

(c) Since $A_{n} \in \mathcal{R}$, $n=1,2,\dots$, such that $A_{n}
\nearrow A$, and $\mu$ is exhaustive on $\mathcal{R}$, then the
sequence $(A_{n})_1^{\infty}$ is $\mu$-Cauchy, i.e.,
$\|\,\mu(A_{m} \triangle A_{n})\| \to 0$ as $n \wedge m \to
\infty$. Considering $m>n$ yields that $A_{m} \triangle A_{n} =
A_{m} \setminus A_{n}$. Thus $\|\,\mu(A_{m} \setminus A_{n})\| \to
0$ as $m \to \infty$. Since $(A_{m} \setminus A_{n}) \nearrow_{m}
(A \setminus A_{n})$, then
$$\|\,\hat{\mu}(A \setminus A_{n})\| = \lim_{m\to \infty} \|\,\mu(A_{m} \setminus A_{n})\|,
\quad \textrm{for every } n\in \mathbb{N},$$ and therefore
$\|\,\hat{\mu}(A \setminus A_{n})\| \to 0$.

(d)  Let $(A_n)_1^\infty$ and $(B_n)_1^\infty$ be two sequences of
sets from $\mathcal{R}_{\sigma}$ and let
$\lim\limits_{n\to\infty}\|\,\hat{\mu}(A_n\triangle B_n)\,\| = 0$.
Then there exist $A_{n,k} \in \mathcal{R}$ and $B_{n,k} \in
\mathcal{R}$, $k=1,2,\dots$, such that $A_{n,k} \nearrow_k A_n$
and $B_{n,k} \nearrow_k B_n$ for each $n\in\mathbb{N}$,
respectively. According to~(\ref{eqhatmu}) for each
$n\in\mathbb{N}$ we have
$$\lim_{k\to\infty} \|\,\mu(A_{n,k})\,\| = \|\,\hat{\mu}(A_n)\,\|
\quad \textrm{and} \quad \lim_{k\to\infty} \|\,\mu(B_{n,k})\,\| =
\|\,\hat{\mu}(B_n)\,\|.$$ Since {\setlength\arraycolsep{2pt}
\begin{eqnarray*} \lim_{n\to\infty}
\lim_{k\to\infty} \|\,\mu(A_{n,k}\triangle B_{n,k})\,\| & = &
\lim_{n\to\infty} \lim_{k\to\infty} \|\,\hat{\mu}(A_{n,k}\triangle
B_{n,k})\,\| \\ & = & \lim_{n\to\infty}
\|\,\hat{\mu}(A_{n}\triangle B_{n})\,\| = 0,
\end{eqnarray*}}then
according to the~(u.s.c.) of $\mu$ on $\mathcal{R}$ (see
Lemma~\ref{lemmaunisubadcont}) we get that for each
$n\in\mathbb{N}$ $$\lim_{k\to\infty} (\|\,\mu(A_{n,k})\,\| -
\|\,\mu(B_{n,k})\,\|) = 0.$$ Then, we have
{\setlength\arraycolsep{2pt}
\begin{eqnarray*}
0 & = & \lim_{n\to\infty} \lim_{k\to\infty} (\|\,\mu(A_{n,k})\,\|
- \|\,\mu(B_{n,k})\,\|) \\ & = & \lim_{n\to\infty}
\left(\lim_{k\to\infty}\|\,\mu(A_{n,k})\,\| -
\lim_{k\to\infty}\|\,\mu(B_{n,k})\,\|\right) \\ & = &
\lim_{n\to\infty} (\|\,\hat{\mu}(A_n)\,\| -
\|\,\hat{\mu}(B_n)\,\|).
\end{eqnarray*}}Thus, according to Lemma~\ref{lemmaunisubadcont} the set function
$\hat{\mu}$ satisfies the~(u.s.c.) on $\mathcal{R}_{\sigma}$.

(e) Let $A_n \in \mathcal{R}_\sigma$, $n=1,2,\dots$, be such that
$A_n \searrow \emptyset$. Then $B_n = A_n \setminus A_{n+1}$,
$n\in \mathbb{N}$, are pairwise disjoint sets from
$\mathcal{R}_\sigma$ and $A_n = \bigcup_{i=n}^{\infty} B_i$. Since
$\hat{\mu}$ is exhaustive on $\mathcal{R}_\sigma$ and has
the~(p.g.p.), then for each $k=2,3,\dots$ there exists an $n_k
> n_{k-1}$, such that
$$\left\|\,\hat{\mu}\left(\bigcup_{i=n_k}^{n_k+p} B_i\right)\right\| < \delta_k \quad \textrm{for each
}\,p=1,2,\dots,$$ 
Thus
$$\left\|\,\hat{\mu}\left(\bigcup_{i=n_j}^{n_{j+1}} B_i\right)\right\| <
\delta_j \quad \textrm{for each }\,j=1,2,\dots,$$ and then
$$\|\,\hat{\mu}(A_{n_k})\| = \left\|\,\hat{\mu}\left(\bigcup_{i=n_k}^{\infty} B_i\right)\right\|
=\left\|\,\hat{\mu}\left(\bigcup_{j=k}^{\infty}\bigcup_{i=n_j}^{n_{j
+1}} B_i\right)\right\| < \delta_{k-1}$$ for each $k=2,3,\dots$.
Since $\delta_k \searrow 0$, then $\|\,\hat{\mu}(A_{n_k})\|\to 0$
as $k\to\infty$. Thus, $\hat{\mu}$ is continuous on
$\mathcal{R}_\sigma$. \qed

Put
$$\mathcal{R}^{*} = \{A; A \subset B \,\textrm{ for some }\, B\in
\mathcal{R}_{\sigma}\}.$$ Obviously, $\mathcal{R}_{\sigma} \subset
\mathcal{R}^{*}$ and $\mathcal{R}^{*}$ is a $\sigma$-ring. For
every $A \in \mathcal{R}^{*}$ define a set function $\mu^{*}:
\mathcal{R}^{*} \to \overline{\Lambda}$ as follows
\begin{equation}\label{mu*}
\mu^{*}(A) = \inf \{\hat{\mu}(B); A \subset B, B \in
\mathcal{R}_{\sigma}\}.
\end{equation} Observe that $\mu^{*} \mid_{\mathcal{R}_{\sigma}} = \hat{\mu}$ and $\mu^{*}$ is monotone. Note
that the $\sigma$-ring $\mathcal{R}^{*}$ is complete with respect
to (Fr\'echet-Nikodym) pseudometric $\rho(A,B) = \mu^{*}(A
\triangle B)$, see~\cite[Corollary 2]{Dobrakov84}. Since
$\hat{\mu}: \mathcal{R}_{\sigma} \to \overline{\Lambda}$ is a
$D_{u}$-submeasure, then clearly $\mu^{*}: \mathcal{R}^{*} \to
\overline{\Lambda}$ satisfies the~(u.s.c.). Note that $\mu^{*}$
need not be necessarily continuous on the whole $\sigma$-ring
$\mathcal{R}^{*}$, but we will show its continuity on
$\mathcal{R}_{0} = \overline{\mathcal{R}}^{\Gamma(\mu^{*})}
\subset \mathcal{R}^{*}$. Also, some other useful properties of
the set function $\mu^{*}$ are summarized in the following lemma.

\begin{lem}\label{lemmamu*}
Let $\mu^{*}$ be defined as in~(\ref{mu*}) and $\mathcal{R}_{0} =
\overline{\mathcal{R}}_{\sigma}^{\Gamma(\mu^{*})}$. Then
\begin{itemize}
\item[(i)] $A \in \mathcal{R}_{0}$ if and only if there exists a
sequence $(A_{n})_1^\infty$ of sets from $\mathcal{R}_{\sigma}$,
such that $\|\,\mu^{*}(A \triangle A_{n})\| \to 0$ as
$n\to\infty$; \item[(ii)] $\mathcal{R}_{0} =
\overline{\mathcal{R}}^{\Gamma(\mu^{*})}$; \item[(iii)] if $A \in
\mathcal{R}_{0}$, then there exists a sequence $(C_{n})_1^\infty$
of sets from $\mathcal{R}_{\sigma}$ with $C_{1} \supset C_{2}
\supset \dots$, such that $A \subset C_{n}$ for every
$n=1,2,\dots$, and $\|\,\mu^{*}(C_{n} \setminus A)\| \to 0$ as
$n\to\infty$; \item[(iv)] $\mu^{*}$ is continuous on
$\mathcal{R}_{0}$.
\end{itemize}
\end{lem}

\proof (i) Let $A \in \mathcal{R}_{0}$ and $\varepsilon >0$.
Suppose that $$\mathcal{V} = \{B; B \in \mathcal{R}^{*},
\|\,\mu^{*}(B)\| \leq \varepsilon\}$$ is an arbitrary neighborhood
of $\emptyset$ in $\mathcal{R}^{*}$. Then the neighborhood $A
\triangle \mathcal{V}$ of $A$ contains an element $E = A \triangle
C \in \mathcal{R}_{\sigma}$, where $C \in \mathcal{V}$. Clearly,
$\|\,\mu^{*}(C)\| \leq \varepsilon$, i.e., $\|\,\mu^{*}(A
\triangle E)\| \leq \varepsilon$.

Now, for a given sequence $(\frac{\varepsilon}{2^{n}})_1^{\infty}$
of positive numbers there exists a sequence $(A_{n})_1^\infty$ of
sets from $\mathcal{R}_{\sigma}$, such that $\|\,\mu^{*}(A
\triangle A_{n})\| \leq \frac{\varepsilon}{2^{n}}$ for
$n=1,2,\dots$. Thus, $\|\,\mu^{*}(A \triangle A_{n})\| \to 0$ as
$n\to\infty$.

Conversely, let $A \in \mathcal{R}^{*}$ and $\|\,\mu^{*}(A
\triangle A_{n})\| \to 0$ as $n\to\infty$ for a sequence
$(A_{n})_1^{\infty}$ of sets from $\mathcal{R}_{\sigma}$. By the
definition of $\mathcal{R}_0$ we have $A \in \mathcal{R}_{0}$.

(ii) Let $\varepsilon > 0$ be chosen arbitrarily and $A \in
\mathcal{R}_{0}$. Then by (i) there exists a sequence
$(A_{n})_1^\infty$ of sets from $\mathcal{R}_{\sigma}$, such that
$\|\,\mu^{*}(A \triangle A_{n})\| \to 0$ as $n\to\infty$.
Accordingly, we may find a positive integer $N$, such that
$\|\,\mu^{*}(A \triangle A_{n})\| < \frac{\varepsilon}{2}$ for
each $n\geq N$. Let $(A_{n,k})_{k=1}^\infty$ be a sequence of sets
from $\mathcal{R}$, such that $A_{n,k} \nearrow_{k} A_{n}$ for
each $n\in \mathbb{N}$. Then by Theorem~\ref{thmhatmu}(c)
$$\lim_{k\to \infty} \|\,\hat{\mu}(A_{n} \triangle A_{n,k})\| =
\lim_{k\to\infty} \|\,\hat{\mu}(A_{n} \setminus A_{n,k})\| = 0,
\quad n=1,2,\dots .$$ Since $\mu^{*} \mid_{\mathcal{R}_{\sigma}} =
\hat{\mu}$, we get
$$\lim_{k\to \infty} \|\,\mu^{*}(A_{n} \triangle A_{n,k})\| = 0,
\quad n=1,2,\dots .$$ As in Theorem~\ref{thmdense} we may prove
that $A \in \overline{\mathcal{R}}^{\Gamma(\mu^{*})}$ and
therefore $\mathcal{R}_{0} \subset
\overline{\mathcal{R}}^{\Gamma(\mu^{*})}$. Also, since
$\mathcal{R} \subset \mathcal{R}_{\sigma}$, then
$\overline{\mathcal{R}}^{\Gamma(\mu^{*})} \subset
\overline{\mathcal{R}}_{\sigma}^{\Gamma(\mu^{*})}$. Hence,
$\mathcal{R}_{0} = \overline{\mathcal{R}}^{\Gamma(\mu^{*})}$. From
Lemma~\ref{lemmadrew1} it follows that $\mathcal{R}_{0}$ is a
ring.

(iii) Since $A \in \mathcal{R}_{0} =
\overline{\mathcal{R}}^{\Gamma(\mu^{*})}$, there exists a sequence
$(A_{n})_1^\infty$ of sets from $\mathcal{R}$, such that
$\|\,\mu^{*}(A \triangle A_{n})\| \to 0$ as $n\to\infty$. Let
$\varepsilon >0$ be arbitrary. From the definition of $\mu^{*}$
and Lemma~\ref{lemma2}(i) it follows that for each $n\in
\mathbb{N}$ there exists a set $F_{n} \in \mathcal{R}_{\sigma}$
such that $A \triangle A_{n} \subset F_{n}$ and
$$\|\,\hat{\mu}(F_{n})\| < \|\,\mu^{*}(A \triangle A_{n})\| +
\frac{\varepsilon}{2^{n}}.$$ Since $\mu^{*}
\mid_{\mathcal{R}_{\sigma}} = \hat{\mu}$, then
\begin{equation}\label{eq1}
\|\,\mu^{*}(F_{n})\| < \|\,\mu^{*}(A \triangle A_{n})\| +
\frac{\varepsilon}{2^{n}},
\end{equation} and we put $G_{n}= \bigcap_{i=1}^{n}(A_{i} \cup F_{i})$.
Clearly, $G_{n} \in \mathcal{R}_{\sigma}$, $n=1,2,\dots$, and
$G_{1} \supset G_{2} \supset \dots$. Also, $$A = (A \setminus
A_{n}) \cup (A \cap A_{n}) \subset (A \setminus A_{n}) \cup A_{n}
\subset A_{n} \cup F_{n}, $$ for each $n\in \mathbb{N}$. Thus, $A
\subset G_{n}$ for each $n\in \mathbb{N}$ and then
$$G_{n}\setminus A \subset (A_{n} \cup F_{n})\setminus A \subset
F_{n}.$$ From monotonicity of $\mu^{*}$ and~(\ref{eq1}) it follows
that $\|\,\mu^{*}(G_{n}\setminus A)\| \to 0$ as $n\to\infty$.

(iv) First we show that $\mu^*$ is exhaustive on $\mathcal{R}_0$.
Suppose the contrary. Since $\mu^*$ has the~(p.g.p.) on
$\mathcal{R}_0$, take the corresponding sequence
$(\delta_k)_1^{\infty}$. Then there exists a positive integer $K$
and a sequence $(A_n)_1^\infty$ of pairwise disjoint sets from
$\mathcal{R}_0$, such that $\|\,\mu^*(A_n)\,\|
> \delta_K$ for each $n\in \mathbb{N}$.
By~(i) for each $n\in \mathbb{N}$ there exists sequence
$(B_{n,l})_{l=1}^\infty$ of sets from $\mathcal{R}_\sigma$, such
that $\|\,\mu^*(A_n\triangle B_{n,l})\,\| \to 0$ for each
$n\in\mathbb{N}$. Thus for each $n\in\mathbb{N}$ there exists a
positive integer $L_n$, such that for each $l\geq L_n$ holds
$\|\,\mu^*(A_n\triangle B_{n,l})\,\| < \delta_{K+3+n}$. Putting
$C_n = B_{n,L_n}$, $n\in \mathbb{N}$ we have $C_n\in
\mathcal{R}_\sigma$ and $\|\,\mu^*(A_n\triangle C_n)\,\| <
\delta_{K+3+n}$ for each $n\in \mathbb{N}$. Since for $n\neq m$
holds $$C_n \cap C_m \subset (A_n \triangle C_n) \cup (A_m
\triangle C_m),$$ then from the~(p.g.p.) $\|\,\mu^*(C_n \cap
C_m)\,\| < \delta_{K+2+n\wedge m}$. Put $$E_1 = C_1, \quad E_n =
\bigcap_{i=1}^{n-1} C_n \setminus C_i, \,\,n\geq 2.$$ Clearly,
$E_n$, $n=1,2,\dots$, are pairwise disjoint sets from
$\mathcal{R}_\sigma$. Since $\mu^*\mid_{\mathcal{R}_\sigma} =
\hat{\mu}$ and $\hat{\mu}$ is exhaustive on $\mathcal{R}_\sigma$,
then there exists a positive integer $N$, such that for each
$n\geq N$ holds $\|\,\mu^*(E_n)\,\| = \|\,\hat{\mu}(E_n)\,\|<
\delta_{K+3}$. Since $$C_n \setminus E_n = \bigcup_{i=1}^{n-1}
(C_i\cap C_n),$$ then for each $n\in \mathbb{N}$ we have
$\|\,\mu^*(C_n\setminus E_n)\,\| < \delta_{K+2}$. Then by~(p.g.p.)
for each $n\geq N$ holds $\|\,\hat{\mu}(C_n)\,\| =
\|\,\mu^*(C_n)\,\| \leq \|\,\mu^*((C_n\setminus E_n)\cup E_n)\,\|
< \delta_{K+1}$. Hence for $n\geq N$ we have the contradiction
$\|\,\mu^*(A_n)\,\| \leq \|\,\mu^*(A_n \triangle C_n)\,\| <
\delta_K$, which proves that $\mu^*$ is exhaustive.

Let $F_n \in \mathcal{R}_0$, $n=1,2,\dots$, be such that $F_n
\searrow \emptyset$. Then $G_n = F_n \setminus F_{n+1}$, $n\in
\mathbb{N}$, are pairwise disjoint sets from $\mathcal{R}_0$, such
that $F_n = \bigcup_{i=n}^{\infty} G_i$. Now in the same way as in
case~(e) of Theorem~\ref{thmhatmu} we obtain that
$\|\,\mu^*(F_n)\|\to 0$ as $n\to\infty$. \qed

Note that $\mu^{*}$ is also order bounded. Now, we are able to
prove the following extension theorem for $D_{u}$-submeasures from
a ring $\mathcal{R}$ to the $\sigma$-ring $\mathcal{R}_{0}$.

\begin{thm}\label{extensionthm}
If $\mu$ is an order bounded exhaustive $D_{u}$-submeasure on a
ring $\mathcal{R}$ of subsets of a set $T\neq \emptyset$, then
there exists a $\sigma$-ring $\mathcal{R}_{0}$ of subsets of $T$,
such that $\mathcal{R} \subset \mathcal{R}_{0}$ and $\mu$ may be
extended to the $D_{u}$-submeasure $\mu^{*}$ on $\mathcal{R}_{0}$,
such that
\begin{itemize}
\item[(a)] $\mathcal{R}_{0} =
\overline{\mathcal{R}}^{\Gamma(\mu^{*})}$; \item[(b)] the
$\sigma$-ring $\mathcal{R}_{0}$ is null-complete with respect to
$\mu^{*}$; \item[(c)] if $\nu$ is a $D_{u}$-submeasure on
$\mathcal{R}_{0}$, such that $\nu \mid_ \mathcal{R} = \mu$, then
for every $A \in \mathcal{R}_{0}$ holds $\|\nu(A)\| =
\|\,\mu^{*}(A)\|$; \item[(d)] the $\sigma$-ring $\mathcal{R}_{0}$
is a null-completion of $\sigma(\mathcal{R})$.
\end{itemize}
\end{thm}

\proof Let $(A_{n})_1^\infty$ be a sequence of sets from
$\mathcal{R}_{0}$, such that $A= \bigcup_{n=1}^{\infty} A_{n}$.
Similarly as in Theorem~\ref{thmdense} we may show that $A \in
\mathcal{R}_{0}=\overline{\mathcal{R}}^{\Gamma(\mu^{*})}$.
Therefore $\mathcal{R}_{0}$ is a $\sigma$-ring containing
$\mathcal{R}$ and $\mu^{*}$ is a $D_{u}$-submeasure on
$\mathcal{R}_{0}$ which is an extension of $\mu$. Thus, the item
(a) is proved.

(b) Let $A \in \mathcal{R}_{0}$ with $\mu^{*}(A)=0$. Then
$\|\,\mu^{*}(A)\|=0$. Since $\mathcal{R}_{0} \subset
\mathcal{R}^{*}$, then $A \in \mathcal{R}^{*}$. Accordingly, $A
\subset C$ for some $C \in \mathcal{R}_{\sigma}$. Then $B \subset
A$ implies $B \subset C \in \mathcal{R}_{\sigma}$. Thus, $B \in
\mathcal{R}^{*}$ and from monotonicity $\|\,\mu^{*}(B)\| \leq
\|\,\mu^{*}(A)\|$ we get $\|\,\mu^{*}(B)\| =0$, and so
$\mu^{*}(B)=0$.

Now we prove that $B \in \mathcal{R}_{0}$. Let $\varepsilon >0$ be
chosen arbitrarily. From the definition of $\mathcal{R}_{0}$ it
follows that there exists $E \in \mathcal{R}$, such that
\begin{equation}\label{eq2}
\|\,\mu^{*}(A \triangle E)\| \leq \varepsilon.
\end{equation} Since $\|\,\mu^{*}(A)\| = \|\,\mu^{*}(B)\| = 0$ and $\mu^{*}$
is monotone, then
\begin{equation}\label{eq3}
\|\,\mu^{*}(A \cup E)\| = \|\,\mu^{*}(A \triangle E)\| =
\|\,\mu^{*}(E)\|,
\end{equation} and
\begin{equation}\label{eq4}
\|\,\mu^{*}(B \cup E)\| = \|\,\mu^{*}(B \triangle E)\| =
\|\,\mu^{*}(E)\|.
\end{equation} Using~(\ref{eq2}), (\ref{eq3}) and~(\ref{eq4}) yields
$$\|\,\mu^{*}(B \triangle E)\| \leq \varepsilon, \quad \textrm{for } E\in \mathcal{R}.$$
Consequently, $B \in \mathcal{R}_{0}$.

(c) Let $\nu$ be a $D_{u}$-submeasure on $\mathcal{R}_{0}$, such
that $\nu \mid_ \mathcal{R} = \mu$ and let $B \in
\mathcal{R}_{\sigma}$. Then there exists a sequence
$(B_{n})_1^\infty$ of sets from $\mathcal{R}$, such that $B_{n}
\nearrow B$. From the definition of $\mu^{*}$ it follows that
$\mu^{*}(B) \leq \nu(B)$. Using~(\ref{eqhatmu}) and
Theorem~\ref{thmRsigma} we may prove that $\mu^{*}(B) = \nu(B)$.
Thus, $\nu \mid_{ \mathcal{R}_{\sigma}} = \hat{\mu}$.

Let $A \in \mathcal{R}_{0}$. Similarly as in
Lemma~\ref{lemmamu*}(iii) there exists a sequence
$(F_{n})_1^\infty$ of sets from $\mathcal{R}_{\sigma}$ with $F_{1}
\supset F_{2} \supset \dots $, such that $A \subset F_{n}$ for
every $n=1,2,\dots$, and
\begin{equation}\label{eq5}
\|\,\mu^{*}(F_{n}\setminus A)\| \to 0 \,\,\textrm{as}\,\,
n\to\infty.
\end{equation} This yields
\begin{equation}\label{eq6}
\|\,\mu^{*}(A)\| = \lim_{n\to \infty} \|\,\hat{\mu}(F_{n})\| =
\lim_{n\to \infty} \|\nu(F_{n})\|.
\end{equation} Let $\varepsilon >0$ be chosen arbitrary.
Since $F_{n}\setminus A \in \mathcal{R}^{*}$, then from the
definition of $\mu^{*}$ it follows that for each $n \in
\mathbb{N}$ there exists $G_{n} \in \mathcal{R}_{\sigma}$, such
that $F_{n} \setminus A \subset G_{n}$ and
$$\|\,\hat{\mu}(G_{n})\| < \|\,\mu^{*}(F_{n} \setminus A)\| + \frac{\varepsilon}{2^{n}}.$$
Consequently, from~(\ref{eq5}) we get $\|\,\hat{\mu}(G_{n})\| \to
0$ as $n\to\infty$. From monotonicity of $\nu$ on $\mathcal{R}$ we
have $\|\nu(F_{n} \setminus A)\| \leq \|\nu(G_{n})\| =
\|\,\hat{\mu}(G_{n})\|$ and therefore $\|\nu(F_{n} \setminus A)\|
\to 0$ as $n\to\infty$. From it follows that $\|\nu(F_{n})\| \to
\|\nu(A)\|$ and from~(\ref{eq6}) we get $\|\nu(A)\| =
\|\,\mu^{*}(A)\|$ for every $A \in \mathcal{R}_{0}$.

(d) Let $A \in \mathcal{R}_{0}$. Then by Lemma~\ref{lemmamu*}(iii)
there exists a sequence $(C_{n})_1^\infty$ of sets from
$\mathcal{R}_{\sigma}$ with $C_{1} \supset C_{2} \supset \dots$,
such that $A \subset C_{n}$ for every $n=1,2,\dots$, and
$\|\,\mu^{*}(C_{n} \setminus A)\| \to 0$ as $n\to\infty$. Let
$C=\bigcap_{n=1}^{\infty} C_{n}$. Then $A \subset C \in
\sigma(\mathcal{R})$ and thus $\|\,\mu^{*}(C \setminus A)\| \leq
\|\,\mu^{*}(C_{n}\setminus A)\|$ for $n=1,2,\dots$. Hence,
$\|\,\mu^{*}(C\setminus A)\| \leq 0$.

Also, $C \setminus A \in \mathcal{R}_{0}$. By
Lemma~\ref{lemmamu*}(iii) there exists a sequence
$(E_{n})_1^\infty$ of sets from $\mathcal{R}_{\sigma}$ with $E_{1}
\supset E_{2} \supset \dots$ and $C\setminus A \subset E_{n}$ for
$n=1,2,\dots$, such that $\|\,\mu^{*}(E_{n}\setminus (C\setminus
A))\|\to 0$ as $n\to \infty$. So,
$$\lim_{n\to \infty}\|\,\mu^{*}(E_{n})\| = \|\,\mu^{*}(C\setminus A)\|
= 0.$$ Now, $$C\setminus A \subset \bigcap_{n=1}^{\infty} E_{n} =
E \in \sigma(\mathcal{R}),$$ and also from monotonicity
$$\|\,\mu^{*}(E)\,\| = \left\|\,\mu^{*}\left(\bigcap_{n=1}^{\infty}
E_{n}\right)\,\right\| \leq \|\,\mu^{*}(E_{n})\,\|, \quad
\textrm{for every }\, n \in \mathbb{N}.$$ From it results that
$\|\,\mu^{*}(E)\| = 0$. Now,
$$C = (C \setminus A) \cup A \subset E \cup A.$$ Since $A \subset
C$, then $A\setminus E \subset C\setminus E$, and since $C \subset
E\cup A$, then $C\setminus E \subset (E\cup A)\setminus E =
A\setminus E$. Thus, $C\setminus E = A \setminus E \subset A
\subset C$ and $C\setminus E$, $E \in \sigma(\mathcal{R})$ and
$$\|\,\mu^{*}(C \setminus (C \setminus E))\| = \|\,\mu^{*}(C \cap E)\|
= 0.$$ Therefore, $\mu^{*}(C \setminus (C \setminus E)) =
\mu^{*}(C \cap E) = 0$, i.e., $\mathcal{R}_{0}$ is a
null-completion of $\sigma(\mathcal{R})$. \qed

\begin{rem}\rm
In Remark~\ref{remark3.6} we have stated that $D$-submeasures do
not satisfy the condition~(a.c.) in general, which seems to play
the crucial role for $\Gamma(\mu)$ to be the $FN$-topology. In
spite of this fact, \textit{is it possible to provide the
(analogous) extension for $D$-submeasures in general}?
\end{rem}

\vspace{5mm}

\noindent \small{Ondrej Hutn\'ik, Institute of Mathematics,
Faculty of Science, Pavol Jozef \v Saf\'arik University in Ko\v
sice, {\it Current address:} Jesenn\'a 5, 040~01 Ko\v sice,
Slovakia,
\newline {\it E-mail address:} ondrej.hutnik@upjs.sk}

\end{document}